\documentclass[a4paper,11pt]{amsart}

\usepackage{lineno,hyperref}

\usepackage{amsfonts,amsmath, amsthm, amscd, amssymb, graphicx, color, mathrsfs, stmaryrd}
\usepackage[utf8]{inputenc}

\def\findemo{\hfill \rule{6pt}{6pt}}
\newcommand{\conv}{\operatorname{conv}}
\newcommand{\spn}{\operatorname{span}}
\newcommand{\dom}{\operatorname{dom}}

\newcommand{\Dz}{\operatorname{Dz}}

\newcommand{\dd}{\operatorname{d}}

\begin{document}

\newtheorem{thm}{Theorem}[section]
\newtheorem{theo}[thm]{Theorem}
\newtheorem{prop}[thm]{Proposition}
\newtheorem{coro}[thm]{Corollary}
\newtheorem{lema}[thm]{Lemma}
\newtheorem{defi}[thm]{Definition}
\newtheorem{ejem}[thm]{Example}
\newtheorem{rema}[thm]{Remark}
\newtheorem{fact}[thm]{Fact}
\newtheorem{open}[thm]{Problem}

\title{Applications of the quantification  \linebreak of super weak compactness}

\author{G. Grelier, M. Raja}
\date{February 28, 2021}

\address{Departamento de Matem\'aticas, Universidad de Murcia, Campus de Espinardo, 30100 Espinardo, Murcia, Spain}

\thanks{The authors were supported by the Grants of Ministerio de Econom\'ia, Industria y Competitividad MTM2017-83262-C2-2-P; and Fundaci\'on S\'eneca Regi\'on de Murcia 20906/PI/18.}

\maketitle

\begin{center}
	{\it Dedicated to our friend Gilles Godefroy with admiration and gratitude.}
\end{center}

\begin{abstract}
We introduce a measure of super weak noncompactness $\Gamma$ defined for bounded linear operators and subsets in Banach spaces that allows to state and prove a characterization of the Banach spaces which are subspaces of a Hilbert generated space. The use of super weak compactness and $\Gamma$ casts light on the structure of 
these Banach spaces and complements the work of Argyros, Fabian, Farmaki, Godefroy, Hájek, Montesinos,\linebreak Troyanski and Zizler on this subject. 
A particular kind of relatively super weakly compact sets, namely uniformly weakly null sets, plays an important role and exhibits connections with Banach-Saks type properties.
\end{abstract}

\section{Introduction}

Along the paper $X$ is a real Banach space, its unit ball is denoted $B_X$ and $X^*$ stands for the dual. In general, our notation is quite standard and the knowledge requirements minimal, however we can address the reader to \cite{banach, JL} for any unexplained notation or concept.
Ultrapowers are a powerful tool to provide brief equivalent definitions of the main notions here (see \cite{heinrich2} for an account of that method in Banach space theory). Here we will consider only ultrafilters on ${\Bbb N}$, although the theory is much richer allowing arbitrary cardinals.
Given a free ultrafilter 
${\mathcal U}$, recall that $X^{\mathcal U}$ is the quotient of $\ell_\infty(X)$ by the subspace of those 
$(x_n)_{n \in {\Bbb N}}$ such that $\lim_{n, \mathcal U} \|x_n\|=0$. A Banach space is said to be {\it super-reflexive} if for some (or, equivalently, any) nontrivial ultrafilter ${\mathcal U}$ on $\Bbb N$, its ultrapower $X^{\mathcal U}$ is reflexive. 
The most representative results on super-reflexive Banach spaces are James' characterizations \cite{James1}, Enflo's uniformly convex renorming \cite{Enflo} and Pisier's applications to Banach valued martingales and renormings with power type modulus \cite{Pisier}.
See the books \cite{beau, banach} for an account on the theory of super-reflexive spaces.\\

Beauzamy \cite{beau1} introduced an operator version of super-reflexivity under the name of {\it uniformly convexifying} property (of an operator), but it was later renamed.
An operator $T: X \rightarrow Y$  is said to {\it super weakly compact} (SWC) if the induced operator $T^{\mathcal U}:
X^{\mathcal U}  \rightarrow Y^{\mathcal U}$  is weakly compact for any ultrafilter ${\mathcal U}$ (equivalently, a free ultrafilter on ${\Bbb N}$). Note that we can think of taking ultrapowers for a fixed ultrafilter ${\mathcal U}$ as a funtor on the category of Banach spaces.
The set of super weakly compact operators is an operator ideal denoted by 
${\mathfrak W}^{super}$. Notably, ${\mathfrak W}^{super}$ is a symmetric ideal, that is, $T \in {\mathfrak W}^{super}$ if and only if  $T^* \in {\mathfrak W}^{super}$.  See  \cite{beau1,beau2, DJP, heinrich} for more properties of 
${\mathfrak W}^{super}$ and its relation with other operator ideals. See also \cite{wenzel} for characterizations in terms of martingale type and cotype, and \cite{CauseyDilworth} for a nonlinear characterization. All the operators considered in this paper are supposed to be linear and bounded.\\

A localized version of super-reflexivity was introduced by the second named author in \cite{raja} for convex sets (and, somehow more generally, for non-linear maps) with the name of finitely dentable sets. The more natural name super weakly compactness was introduced in \cite{Cheng}. Given a set $A \subset X$ we will denote $A^{\mathcal U}$ the subset of $X^{\mathcal U}$ whose elements have a representative in $A^{\Bbb N}$. A set $A \subset X$ is said to be {\it relatively super weakly compact} (relatively SWC) if $A^{\mathcal U}$ is a relatively weakly compact subset of $X^{\mathcal U}$ for some (or, equivalently, any) free ultrafilter ${\mathcal U}$.  Moreover, $A \subset X$ is said to be {\it super weakly compact} (SWC, of course) if it is relatively super weakly compact and weakly closed. The class of SWC sets lies strictly between the norm compact and the weakly compact subsets.
The theory of SWC sets has been developed during the last 15 years in a series of papers \cite{raja, Cheng, Cheng2, raja2, YLC, Cheng3, KT, KT2, LR, GR}.\\

Super weak compactness is more widespread than it may appear. For instance, any weakly compact operator with range $L_1(\mu)$ ($\mu$ any measure) or domain $C(K)$ ($K$ any Hausdorff compact) is super weakly compact, see \cite[Proposition 6.1]{LR}. Actually, some results in Banach space theory could be understood in terms of super weak compactness. As for instance, the classic Szlenk result establishing that a weakly convergent sequence in $L_1(\mu)$ has a subsequence whose {\it Cesàro means} converge (to the same limit) is a consequence of two facts: the weakly compact subsets of $L_1(\mu)$ are SWC; and the  SWC sets have the {\it Banach-Saks property} \cite[Corollary 6.3]{LR}. In relation with the structure of nonseparable Banach spaces, the second named author \cite{raja2} showed the role of SWC generated (and strongly generated) Banach spaces among the subspaces of Hilbert generated spaces.\\

The aim of this paper is to show that, actually, super weak compactness and, particularly, its quantification, may cast light on the structure of the subspaces of Hilbert generated Banach spaces. Indeed, we have realized that several ``technical hypotheses'' in papers of Troyanski \cite{tro2}, Argyros and Farmaki \cite{farma}, and the series by Fabian, Godefroy, Hájek, Montesinos and Zizler \cite{FHZ, FGZ, FGHZ, FHMZ2}
on the structure of Hilbert generated spaces and uniformly Gâteaux renorming,
can be understood in terms of a quantified version of super weak compactness.\\

For a better understanding of our main result, we will state firstly the ``non uniform version'' with the help of a measure of  weak noncompactness. Let $A \subset X$ be a bounded set, then take
$$ \gamma(A) = \inf\{ \varepsilon>0: \overline{A}^{w^*} \subset X + \varepsilon B_{X^{**}}\} .$$
We have that a set $A$ is relatively weakly compact if and only if $\gamma(A)=0$. This measure has been studied in \cite{FHMZ, granero, CMR}, see also \cite[Section 3.6]{BIO}, and there are several measures of weak noncompactness that turn out to be equivalent \cite{AC}.

\begin{theo}[\cite{BRW,FMZ}]\label{WCG}
For a Banach space $X$ the following statements are equivalent:
\begin{itemize}
\item[(i)] $X$ is a subspace of a WCG space;
\item[(ii)]  $(B_{X^*},w^*)$ is an Eberlein compact;
\item[(iii)] For every $\varepsilon>0$ there are sets $(A_n^\varepsilon)$ such that $B_X=\bigcup_{n=1}^\infty A_n^\varepsilon$ and $\gamma(A_n^\varepsilon) < \varepsilon$.
\end{itemize}
\end{theo}

The equivalence (i)$\Leftrightarrow$(ii) is due to Benyamini, Rudin and Wage \cite{BRW}. The inner characterization (iii) was obtained by Fabian, Montesinos and Zizler \cite{FMZ}. Recall that WCG stands for {\it weakly compactly generated}, that is, a Banach space that contains a weakly compact subset whose linear span is dense. Thanks to the celebrated interpolation result of Davis, Figiel, Johnson and Pełczyński \cite{DJLP} (see also \cite[Theorem 13.22]{banach}), a Banach space $X$ is WCG if and only if there exists a reflexive space $Z$ and an operator $T:Z \rightarrow X$ with dense range. Moreover, if the space $Z$ can be taken a Hilbert space, we say that $X$ is {\it Hilbert generated}. The name 
{\it Eberlein} applies to the compact spaces which are homeomorphic to a weakly compact set of a Banach space. It is well known after Amir and Lindenstrauss (see \cite[Corollary 13.17]{banach}, for instance) that an Eberlein compact embeds as a weakly (equivalent, bounded and pointwise) compact subset of $c_0(I)$ for $I$ large enough. If such an embedding can be done into a Hilbert space $\ell_2(I)$, then the compact is said to be {\it uniformly Eberlein}.
Note that the third statement in Theorem \ref{WCG} is actually an internal characterization as it is written in terms of the space $X$, not an over-space or its dual. A different matter is if the computation of $\gamma$ requieres the use of the over-space $X^{**}$, as we will see later there are equivalent definitions of $\gamma$ that does not appeal to the bidual space.\\

Let us prepare the way to state the uniform analogue of Theorem  \ref{WCG}. We will requiere the following {\it measure of super weak noncompactness}: for a bounded set $A \subset X$ take
$$ \Gamma(A) := \gamma(A^{\mathcal U}) $$
where ${\mathcal U}$ is a free ultrafilter and $\gamma$ is computed in $X^{\mathcal U}$. Later we will see that 
$\Gamma$ does not depend, essentially, on the choice of the ultrafilter ${\mathcal U}$. Obviously, we have that $A$ is relatively SWC if and only if $\Gamma(A)=0$, and an operator $T:X \rightarrow Y$ is SWC if and only if $\Gamma(T(B_X))=0$. Now we are ready to state our main result. Please, note the parallelism  with Theorem~\ref{WCG}.

\begin{theo}\label{main}
For a Banach space $X$ the following statements are equivalent:
\begin{itemize}
\item[(i)] $X$ is a subspace of a Hilbert generated space;
\item[(ii)] $(B_{X^*},w^*)$ is a uniform Eberlein compact;
\item[(iii)] For every $\varepsilon>0$ there are sets $(A_n^\varepsilon)$ such that $B_X=\bigcup_{n=1}^\infty A_n^\varepsilon$ and $\Gamma(A_n^\varepsilon) < \varepsilon$.
\end{itemize}
\end{theo}

Again, the equivalence (i)$\Leftrightarrow$(ii) goes back to Benyamini, Rudin and Wage \cite{BRW}. Moreover, there are also some other equivalent conditions in terms of uniformly Gâteaux renorming coming from \cite{FHZ, FGZ, FGHZ, FHMZ2} (see also \cite{BIO}) that would spoil the nice analogy with Theorem \ref{WCG}. We will say more on this matter later.
Also, in relation with our Theorem \ref{main}, we will prove that we can change $B_X$ in (iii) by any linearly dense subset of $X$, see Theorem \ref{main2}. Also, in order to apply statement (iii) is quite relevant the fact that $\Gamma$ can be computed in several fashions, some of them without ultrapowers neither over-spaces, see Proposition \ref{superweaknoncom}.\\

Previous works on uniformly Gâteaux renorming by Fabian, Godefroy, Hájek and Zizler \cite{FGHZ}, as well as early results by Troyanski \cite{tro2}, unawarely contain estimations of $\Gamma$. The explanation will come through the following result.

\begin{prop}\label{Banach-Saks}
Let $A \subset X$ a bounded subset and consider the two following numbers:
\begin{itemize}
\item[$(\varepsilon_1)$] is the infimum of the $\varepsilon>0$ such that there is $n_1 \in {\Bbb N}$ such that for every $x^* \in B_{X^*}$ then
$$ | \{ x \in A: |x^*(x)| > \varepsilon \} | \leq n_1 ;$$
\item[$(\varepsilon_2)$]  is the infimum of the $\varepsilon>0$ such that there is $n_2 \in {\Bbb N}$ such that for any finite set $B \subset A$ with $|B| \geq n_2$ then
$$ \left\| \,  \frac{1}{|B|} \sum_{x \in B} x \,\, \right\| < \varepsilon. $$
\end{itemize}
Then $\varepsilon_1=\varepsilon_2$ and in such a case $\Gamma(A) \leq \varepsilon_1$.
\end{prop}

The sets satisfying the statements of Proposition \ref{Banach-Saks} with $\varepsilon_1=\varepsilon_2=0$ will be called {\it uniformly weakly null} sets. Note that a uniformly weakly null set becomes SWC by adding $\{0\}$. Together with unit balls of super-reflexive spaces, uniformly weakly null sets are the most prototypical examples of SWC sets. As we will see, SWC sets with some reasonable discreteness assumption are uniformly weakly null. Note that the second statement ($\varepsilon_2$) is a sort of uniform Banach-Saks property (with unique limit $0$). That will allow us to apply results of infinite combinatorics, such as the Erdös-Magidor \cite{ErMa} and Mercourakis \cite{merco} selections.\\

The structure of the paper is as follows. Section 2 covers some new aspects on super weak compactness with special emphasis on uniformly weakly null sets. There, we show that the Eberlein-\v{S}mulian theorem fails for super weak compactness. In section 3, we study the properties of the measure $\Gamma$ and we show several different ways to compute or estimate it. Our results depends on some equivalent forms for $\gamma$, that may be of independent interest. Section 4 deals with the application of $\Gamma$ to operators. In particular, we prove a quantified version of the symmetry of the 
ideal ${\mathfrak W}^{super}$, as well as a quantified version 
of Beauzamy's renorming to make uniformly convex a super weakly compact operator. Section 5 is devoted to the proof of Theorem \ref{main}, actually a more general version, and an application to Jordan algebras. Finally, section 6 contains more results on uniformly weakly null sets. We investigate when a Schauder basis is uniformly weakly null, and the relation of uniformly weakly null sets with the representation of uniformly Eberlein compacts. Let us point out that the atypically long list of references is a consequence of the, hitherto unnoticed, transversality of super weak compactness in Banach space theory.

\section{Some remarks on super weak compactness}

One big issue dealing with super weak compactness is the lack of representation for $(X^{\mathcal U})^*$, except when $X$ is super-reflexive. In such a case, $(X^{\mathcal U})^*=(X^*)^{\mathcal U}$. Otherwise, 
$(X^*)^{\mathcal U}$ is a proper subspace of  $(X^{\mathcal U})^*$. The identification is done by the 
assignment 
$$ \langle (x_n^*) , (x_n) \rangle := \lim_{n ,{\mathcal U}} x_n^*(x_{n}) . $$
A subset $B \subset B_{X^*}$ is called a {\it boundary} (of $X$) if for every $x \in X$ there exists $x^* \in B$ such that
$\|x\|=x^*(x)$.
A norming subspace $Z \in X^*$ is called a boundary if $B_Z= Z \cap B_{X^*}$ is a boundary in the previous sense. We have the following.

\begin{theo}
	Let $X$ be  Banach space and let ${\mathcal U}$ be  any free ultrafilter. Then $(X^*)^{\mathcal U}$ is a boundary for $X^{\mathcal U}$. Therefore, the relatively weakly compact subsets of  $X^{\mathcal U}$ are exactly those which are relative compact for the topology of pointwise convergence on the elements from  $(X^*)^{\mathcal U}$.
\end{theo}

\noindent
{\bf Proof.} 
 Indeed, for every $n \in {\Bbb N}$ take $x^*_n$ such that $x_n^*(x_n)=\|x_n\|$. The second statement comes from Pfitzner's solution \cite{Pfitzner} to Godefroy's boundary problem.\findemo\\

The most typical example of SWC set is the unit ball of a super-reflexive Banach space. Now we will introduce another family of (relatively) SWC sets which will be more relevant to the results of this paper. 
We say that a subset $A \subset X$ is {\it uniformly weakly null} if for every $\varepsilon >0$ there is $n(\varepsilon) \in {\Bbb N}$ such that, for every $x^* \in B_{X^*}$,
$$ | \{ x \in A : |x^*(x)| > \varepsilon \} |  \leq n(\varepsilon) .$$
Note that any sequence made of different points of a uniformly weakly null set is a weakly null sequence. Therefore, uniformly weakly null sets are relatively weakly compact (and become weakly compact just by adding $0$). We have something better.

\begin{theo}\label{supernull}
	Let $A \subset X$ be a uniformly weakly null set and let ${\mathcal U}$ be  any free ultrafilter. Then $A^{\mathcal U}$ is uniformly weakly null in $X^{\mathcal U}$ and, therefore, $A$ is relatively super weakly compact in $X$.
\end{theo}

\noindent
{\bf Proof.} Let $\overline{x}_1, \dots, \overline{x}_n \in A^{\mathcal U}$ be different vectors,  $\overline{x}^* \in B_{(X^{\mathcal U})^*}$ and  $\varepsilon>0$ such that $|\overline{x}^*(x_k)|>\varepsilon$ for every $1 \leq k \leq n$. We claim that for 
$\varepsilon'<\varepsilon$, there are different elements $x_1, \dots, x_n \in A$ and $x^* \in B_{X^*}$ with $|x^*(x_k)|>\varepsilon'$.
Indeed, the proof that $X$ is finitely representable in $X^{\mathcal U}$ (see \cite[p. 222]{beau} for instance), provides those 
$x_1, \dots, x_n \in X$ in such a way that $\overline{Y} =\spn\{ \overline{x}_1, \dots, \overline{x}_n\}$ and 
$Y = \spn\{ x_1, \dots, x_n \}$ are $\varepsilon/\varepsilon'$-isomorphic. 
Moreover, the vector $x_k$ is found on the ``coordinates'' of $\overline{x}_k$, so we may assume $x_k \in A$ for all $k$.
Then $T: Y \rightarrow \overline{Y}$ be the isomorphism. Let $x^*$ be the Hahn-Banach extension of $(\varepsilon'/\varepsilon) \overline{x}^*  \circ T$. Then,
$x^* \in B_{X^*}$ and $|x^*(x_k)|>\varepsilon'$ for all $1 \leq k \leq n$ as desired.
That claim shows that $A^{\mathcal U}$ have to be uniformly weakly null.
Now we have $A^{\mathcal U}$ is weakly compact in $X^{\mathcal U}$ and thus $A$ is SWC.\findemo\\

A sequence $(x_n)$ that is a uniformly weakly null set is called {\it uniformly weakly null sequence}.
A sequence $(x_n)$ is {\it uniformly weakly convergent} to $x$ if $(x_n-x)$ is a uniformly weakly null sequence. The fact that a uniformly weakly convergent sequence together its limit is a super weakly compact set was noted in \cite{Cheng3}. Uniformly weakly convergent sequences are closely related to the Banach-Saks property. 
A sequence $(x_n)$ is said to be {\it Cesàro convergent} if the sequence of its arithmetic means 
$$ n^{-1} \sum_{k=1}^n x_k $$
converges (in norm) to some $x \in X$. A set $A \subset X$ is said to have the {\it Banach-Saks property} if every sequence $(x_n) \subset A$ has a Cesàro convergent subsequence. Recall that relatively SWC sets are Banach-Saks \cite[Corollary 2.4]{LR}. The relations between both properties is an interesting topic, although  
we will not deal here with the Banach-Saks property in general. Let us introduce the following ``ephemeral'' definition. A set $A \subset X$ is said to be {\it uniformly Banach-Saks null} if for every $\varepsilon>0$ there is $n(\varepsilon)$ such that whenever $B \subset A$ is finite with $|B| \geq n(\varepsilon)$ then
$$ |B|^{-1} \left\| \sum_{x \in B} x \right\| < \varepsilon .$$
Proposition \ref{Banach-Saks} has the following consequence.

\begin{coro}
	Let $A \subset X$ be a bounded subset. Then $A$ is uniformly Banach-Saks null if and only if it is uniformly weakly null.
\end{coro}

\noindent
{\bf Proof of Proposition \ref{Banach-Saks}.} Let $r>0$ such that $A \subset r B_X$.
Take $\varepsilon>\varepsilon_1$ and fix the corresponding number $n_1$. For $n > n_1$ and any $B \subset A$ with $|B|=n$ we have
$$ |x^*( \sum_{x \in B} x) | < n_1 r + (n-n_1) \varepsilon $$ 
for every $x^* \in B_{X^*}$. Therefore
$$ n^{-1} \left\| \sum_{x \in B} x \right\| < \frac{n_1 r}{n} + \left(1- \frac{n_1}{n} \right) \varepsilon $$ 
Since the bound can be taken arbitrarily closed to $\varepsilon$ independently from $B$ if $n$ is large enough, we have that $\varepsilon_2 \leq \varepsilon_1$. That proves the equality  $\varepsilon_1 = \varepsilon_2$ in case $\varepsilon_1 =0$. Assume now that $\varepsilon_1 >0$ and take $0< \varepsilon < \varepsilon_1$. Then, for every $n \in {\Bbb N}$ we can find $C \subset A$ with $|C|=2n$ and $x^* \in B_{X^*}$ such that $x^*(x) > \varepsilon$ or $x^*(x) < -\varepsilon$ for all $x \in C$. Since at least one half of the elements satisfies the same inequality, we may find $B \subset C$ such that $|B|=n$ and
$$ |x^*( \sum_{x \in B} x )| > n \varepsilon. $$
Therefore, we have
$$ n^{-1} \left\| \sum_{x \in B} x \right\| > \varepsilon, $$
that implies $\varepsilon_2 \geq \varepsilon_1$. Now, note that the first statement implies 
$$ \overline{A}^{w^*} \subset A \cup \varepsilon_1 B_{X^{**}} \subset X + \varepsilon_1 B_{X^{**}}  $$
and so $\gamma(A) \leq \varepsilon_1$. In order to pass to $\Gamma$, just follow the ideas in the proof of Theorem \ref{supernull} or check that the property of the second statement is stable by ultrapowers. In any case, we get that 
$\Gamma(A) \leq \varepsilon_1$.\findemo

\begin{rema}\label{norming}
The proof of the equivalence shows that it is enough to check condition ($\varepsilon_2$) for $x^*$ from a norming subset of $B_{X^*}$.
\end{rema}

Mercourakis \cite{merco} improvement of the Erdös-Magidor \cite{ErMa} dichotomy for bounded sequences can be stated in this way (see also \cite{saks} for related results and references).

\begin{theo}[Mercourakis]\label{Mercou}
	Let $(x_n) \subset X$ a bounded sequence. Then there exists a subsequence $(x_{n_k})$ of $(x_n)$ for which one of the following statements holds:
	\begin{itemize}
		\item[(i)] either, $(x_{n_k})$ is uniformly weakly convergent;
		\item[(ii)] or, no subsequence of $(x_{n_k})$ is Cesàro convergent. 
	\end{itemize}
\end{theo}

The celebrated Eberlein-\v{S}mulian theorem, see \cite{banach} for instance, says that weak compactness is determined by sequences.
As an application, we get that there is no Eberlein-\v{S}mulian for super weak compactness. That is, the fact that every sequence has a relatively SWC subsequence does not imply that the set is relatively SWC.

\begin{coro}\label{Mercou2}
	Let $A \subset X$ be a relatively super weakly compact set. Then every sequence $(x_n) \subset A$ contains a uniformly weakly convergent subsequence. However, this property does not characterize the super weakly compactness. Actually, it characterizes the Banach-Saks property.  
\end{coro}

\noindent
{\bf Proof.} For a Banach-Saks set the dichotomy \ref{Mercou} always produces a uniformly weakly convergent subsequence. On the other hand, every uniformly convergent sequence is Cesàro convergent. Therefore, the Banach-Saks property is  characterized by sequences. The other statements follow from the fact that relatively SWC sets are Banach-Saks and there exist Banach-Saks sets which are not relatively SWC \cite[Corollary 2.5]{LR}.\findemo

\section{Different ways to quantify SWC}

Measures of noncompactness can be defined in very general settings. Here we will restrict ourselves to the frame of topological vector spaces. Let $X$ be a topological vector space and let ${\mathfrak K}$ be a vector bornology of compact subsets (that just means the class is stable under some elementary operations). A {\it measure of noncompactness} associated to ${\mathfrak K}$ is a nonnegative function $\mu$ defined on the bounded subsets of $X$ that satisfies the following properties:
\begin{enumerate}
\item $\mu(\overline{A})=\mu(A)$ 
\item $\mu(A)=0$ if and only if $\overline{A} \in {\mathfrak K}$
\item $\mu(A \cup B) = \max\{ \mu(A), \mu(B) \}$
\item $\mu(\lambda A) = |\lambda| \mu(A)$ for all $\lambda \in {\Bbb R}$
\item $\mu(A+B) \leq \mu(A) + \mu(B)$
\item there exists $k>0$ such that $\mu(\conv(A)) \leq k \, \mu(A)$
\end{enumerate}

This list of conditions comes from the usual requirements in literature \cite{ADL} and some properties enjoyed by several measures that are interesting for Banach space geometry, such as $\gamma$ or the family of measures introduced in \cite{LPR} in relation with the Szlenk index.
Condition (6) is usually the most tricky and necessarily requieres that the class ${\mathfrak K}$ be stable by closed convex hulls (Krein-type theorem).\\

The quantification of the super weak non-compactness is linked to the quantification of the weak non-compactness. 
De Blasi (see \cite{AC}, for instance) introduced a measure of weak noncompactness $\omega$ as follows 
$$ \omega(A) = \inf\{ \varepsilon>0: \exists K \subset X \mbox{ weakly compact}, A \subset K + \varepsilon B_X \} .$$
It is not hard to check that $\omega$ enjoys all the properties above. In particular, we have 
$$\omega(\conv(A)) = \omega(A),$$ 
that is, its ``convexifiability constant'' is $1$.
Another quite natural way to measure weak noncompactness, is the function 
$\gamma$ mentioned in the introduction
$$ \gamma(A) =  \inf\{ \varepsilon>0: \overline{A}^{w^*} \subset X + \varepsilon B_{X^{**}}\}  =
\sup \{ d(X,x^{**}): x^{**} \in \overline{A}^{w^*} \} .$$
It is easy to check that $\gamma(A) \leq \omega(A)$ for any bounded set $A \subset X$. However, there is no constant $c>0$ such that $\omega(A) \leq c \, \gamma(A)$ in general, see \cite[Corollary 3.4]{AC}. That fact says that $\omega$ and $\gamma$ are not {\it equivalent}.
The measure $\gamma$ was introduced in \cite{FHMZ} where the authors also proved (\cite{granero} independently, see also \cite[Theorem 3.64]{BIO}) that 
$$\gamma(\conv(A)) \leq 2 \, \gamma(A)$$  
for any bounded $A \subset X$. Notably, there are many different equivalent ways to deal with $\gamma$ which are interesting to us because they have a ``super'' version. The following contains the quantified version of two classic James' characterizations of relative weak compactness together with the quantified version of Grothendieck's commutation of limits criterion.

\begin{prop}\label{weaknoncom}
Let $A \subset X$ be a bounded set. Consider the following numbers:
\begin{itemize}
\item[$(\gamma_1)$] $= \gamma(A)$;
\item[$(\gamma_2)$] the supremum of the numbers $\varepsilon \geq 0$ such that there are sequences $(x_n) \subset A$ and $(x_n^*) \subset B_{X^*}$ such that $x^*_n(x_m)=0$ if $m < n$ and  $x^*_n(x_m) \geq \varepsilon$ if $m \geq n$;
\item[($\gamma_3$)]  the supremum of the numbers $\varepsilon>0$ such that for any $n \in {\Bbb N}$ there are $x_1,\dots,x_n \in C$ such that $\dd( \conv\{x_1,\dots,x_k\},  \conv\{x_{k+1},\dots,x_n\}) \geq \varepsilon$ for all $k=1,\dots, n-1$; 
\item[$(\gamma_4)$] the infimum of the numbers $\varepsilon \geq 0$ such that 
$$ |\lim_n \lim_m x^*_n(x_m) - \lim_m \lim_n x^*_n(x_m) | \leq \varepsilon $$
whenever $(x_n) \subset A$, $(x_n^*) \subset B_{X^*}$ and the iterated limits exist.
\end{itemize}
Then $\gamma_1 \leq \gamma_2 \leq \gamma_3 \leq \gamma_4 \leq 2 \gamma_1$.
\end{prop}

\noindent
{\bf Proof.} Take $\varepsilon < \gamma(A)$ and let $x^{**} \in \overline{A}^{w^*}$ with $d(X,x^{**})>\varepsilon$. We will build sequences satisfying the second statement for such an $\varepsilon$. Indeed, there exists $x^*_1  \in B_{X^*}$ with $|x^{**}(x^*_1)|>\varepsilon$. Now take $x_1 \in A$ such that $x^*_1(x_1) \geq \varepsilon$. 
Assume we have $x_k$ and $x_k^*$ already built for $1 \leq k < n$ and it is satisfied $x^{**}(x_k^*) > \varepsilon$.
An elementary application of Helly's theorem \cite[p. 159]{banach} to $X^{**}$ allows us to choose $x_n^* \in B_{X^*}$ such that $x_n^*(x_k) =0$ for $1 \leq k <n$ and $x^{**}(x_n^*) > \varepsilon$. Now we take
$$ x_n \in A \cap \{ x \in X: x^{*}_k(x) > \varepsilon, 1 \leq k \leq n \} $$
since the set is nonempty. That finishes the construction of the sequence and proves 
$\gamma_1 \leq \gamma_2$.
The inequality $\gamma_2 \leq \gamma_3$ follows straight. In order to prove $\gamma_3 \leq \gamma_4$, 
take $\varepsilon < \gamma_3$, and sequence $(x_n)$ as in the statement ($\gamma_3$).
For every $n \in {\Bbb N}$, take $x^*_n \in B_{X^*}$ such that 
$x^*_n(y) \leq \varepsilon + x^*_n(z)$ for every $y \in \conv\{x_1,\dots, x_n\}$ and $z \in \conv\{x_{n+1}, x_{n+2}, \dots\}$.
The sequences satisfies the following property
$$ x_n^*(x_p) \leq \varepsilon + x_n^*(x_q) $$
whenever $p \leq n < q$. Passing to a subsequence, we may assume the existence of the limits $\lim_n x^*_n(x_m)$ and $\lim_m x^*_n(x_m)$, as well as the existence of the iterated limits. In such a case we will get
$$  \lim_m \lim_n x^*_n(x_m) \leq \varepsilon +  \lim_n \lim_m x^*_n(x_m)  $$
which implies $\varepsilon \leq \gamma_4$, and therefore $\gamma_3 \leq \gamma_4$.
Finally, $\gamma_4 \leq 2 \gamma_1$ is proved in \cite{AC}.\findemo\\

Now we will state the ``super'' version of Proposition \ref{weaknoncom}, for which we prefer to avoid a uniform version of {\it Grothendieck's commutation of limits} (fourth statement).

\begin{prop}\label{superweaknoncom}
Let $A \subset X$. Consider the following numbers:
\begin{itemize}
	\item[($\Gamma_1$)] $= \gamma(A^{\mathcal U})$ measured in $X^{\mathcal U}$ for ${\mathcal U}$ a free ultrafilter;
	\item[($\Gamma_2$)] the infimum of the numbers $\varepsilon>0$ such that there are no arbitrarily long sequences $(x_k)_1^n \subset A$,  $(x_k^*)_1^n \subset B_{X^*}$ with $x^*_k(x_j)=0$ if $j < k$ and $x^*_k(x_j)> \varepsilon$ if $j \geq k$;
	\item[($\Gamma_3$)]  the supremum of the numbers $\varepsilon>0$ such that for any $n \in {\Bbb N}$ there are $x_1,\dots,x_n \in C$ such that $\dd( \conv\{x_1,\dots,x_k\},  \conv\{x_{k+1},\dots,x_n\}) \geq \varepsilon$ for all $k=1,\dots, n-1$; 
\end{itemize} 
Then $\Gamma_1 \leq \Gamma_2 \leq \Gamma_3 \leq 2\Gamma_1$.
\end{prop}

\noindent
{\bf Proof.} The fact $\Gamma_1 \leq \Gamma_2$ follows straight by applying finite representatibity to inequality 
$\varepsilon_1 \leq \varepsilon_2$ in Proposition \ref{weaknoncom}. It is quite easy to get $\Gamma_2 \leq \Gamma_3$, and $ \Gamma_3 \leq 2\Gamma_1$ follows using the standard ultrapower technique, (see also Theorem \ref{quantify-swc} below where the convex case is considered).
\findemo\\

Recall that $\Gamma_1$ is the measure introduced at the introduction
$$ \Gamma(A) := \gamma(A^{\mathcal U}) $$
that depends on the choice of ${\mathcal U}$. From now on, we will assume the free ultrafilter ${\mathcal U}$ is fixed  when speaking of $\Gamma$ or dealing with the ultrapowers.
Note that the equivalent measures $\Gamma_2$ and $\Gamma_3$ does not depend on any ultrafilter. 
Moreover, $\Gamma_3$ does not involves explicitly the dual space. In next section we will use 
$\Gamma_s(A)=\Gamma_2$ as an alternative to $\Gamma(A)$.

\begin{prop}
Let $T: X \rightarrow Y$ be an operator and let $A \subset X$ be a bounded set. Then $\Gamma(T(A)) \leq  \|T\| \,\Gamma(A)$.
\end{prop}

\noindent
{\bf Proof.} Firstly, we will prove a similar statement for $\gamma$.
Consider $T^{**}: X^{**} \rightarrow Y^{**}$ which is weak$^*$ to weak$^*$ continuous. For any bounded set $A \subset X$ we have
$$ \overline{T(A)}^{w^*} = T^{**}(\overline{A}^{w^*}) \subset T(X) + \varepsilon T^{**}(B_{X^{**}}) \subset Y + \varepsilon \|T\| B_{Y^{**}}$$
where $\varepsilon>\gamma(A)$. Therefore $\gamma(T(A)) \leq \|T\| \, \gamma(A)$. In order to prove the statement for $\Gamma$, consider the induced operator $T^{\mathcal U}: X^{\mathcal U} \rightarrow Y^{\mathcal U}$. Then we have
$$ \Gamma(T(A)) = \gamma(T^{\mathcal U}(A^{\mathcal U})) \leq \|T^{\mathcal U} \, \| \gamma(A^{\mathcal U}) = \|T\| \, \Gamma(A) ,$$
as we wished.\findemo\\

In order to state the results from our paper \cite{GR} that we will need later, it is necessary to introduce a certain number of quantities related to sets in Banach spaces. 
Let us denote by ${\Bbb H}$ the set of all the {\it open half-spaces} of $X$, that is, all the sets of the form
$H=\{x\in X,\ x^*(x)>\alpha\}$, with $x^*\in X^*$ and $\alpha \in {\Bbb R}$.
A {\it slice of} of $D \subset X$ is a set of the form $D \cap H \not = \emptyset$, where $H\in {\Bbb H}$.
We say that a bounded closed convex set $C \subset X$ is {\it dentable} if for any nonempty closed convex subset $D \subset C$ has (nonempty) slices of arbitrarily small diameter.
If $C$ is dentable we may consider the following {\it set derivation}:
\[ [D]'_{\varepsilon} = \{ x \in D: \mbox{diam}(D \cap H)>\varepsilon,\ \text{for any}\ H \in {\Bbb H}\ \text{s.t.}\
x \in H \} .\]
Clearly, $[D]'_{\varepsilon}$
is what remains of $D$ after removing all the slices of $D$ of diameter at most
$\varepsilon$. Consider the sequence of sets defined by $[C]_{\varepsilon}^{0}=C$ and, for every $n \in {\Bbb N}$, inductively by
\[ [C]_{\varepsilon}^{n}=[[C]_{\varepsilon}^{n-1}]'_{\varepsilon}. \]
If there is an $n$ in ${\Bbb N}$ such that $[C]_{\varepsilon}^{n-1} \not = \emptyset$ and  $[C]_{\varepsilon}^{n}  = \emptyset$ we set $Dz(C,\varepsilon)=n$. We say that $C$ is {\it finitely dentable} if $Dz(C,\varepsilon)$ is finite for every $\varepsilon>0$.
Given a convex set $C \subset X$, let us denote by $\mbox{Dent}(C)$ the infimum of the numbers $\varepsilon>0$ such that $C$ has nonempty slices contained in balls of radius less than $\varepsilon$, and take 
$\Delta(C)= \sup\{ \mbox{Dent}(B): B \subset C\}$. The measure $\Delta$ was introduced in \cite{CPR} as a way to quantify the lack of {\it Radon-Nikodym property} (RNP).
Let $\varepsilon>0$.
A function $f: X \rightarrow \overline{{\Bbb R}}$ is said to be $\varepsilon$-{\it uniformly convex} with respect to some metric $d$ if 
there is $\delta>0$ such that whenever $d(x,y) \geq \varepsilon$, then  
$$f \left( \frac{x+y}{2} \right) \leq  \frac{f(x)+f(y)}{2}  - \delta. $$
No mention to an explicit metric $d$ means that we are using the norm metric.
The function is said to be just {\it uniformly convex} if it is $\varepsilon$-uniformly convex for all $\varepsilon>0$.

\begin{theo}[\cite{GR}]\label{quantify-swc}
	Let $C \subset X$ be a bounded closed convex subset. Consider the following numbers:
	\begin{itemize}	
		\item[$(\eta_1)$] $=\Gamma(C)$;
		\item[$(\eta_2)$] the supremum of the numbers $\varepsilon>0$ such that for any $n \in {\Bbb N}$ there are $x_1,\dots,x_n \in C$ such that $\dd( \conv\{x_1,\dots,x_k\},  \conv\{x_{k+1},\dots,x_n\}) \geq \varepsilon$ for all $k=1,\dots, n-1$;
		\item[$(\eta_3)$] the supremum of the $\varepsilon>0$ such that there are $\varepsilon$-separated dyadic trees in $C$ of arbitrary height;
		\item[$(\eta_4)$] $=\Delta(C^{\mathcal U})$;
		\item[$(\eta_5)$] the infimum of the $\varepsilon>0$ such that $\Dz(C,\varepsilon) <\omega$;
		\item[$(\eta_6)$] the infimum of the $\varepsilon>0$ such that $C$ supports a convex bounded $\varepsilon$-uniformly convex function.
	\end{itemize}
	Then $\eta_1 \leq \eta_2 \leq 2\eta_3 \leq 2\eta_4 \leq 2\eta_1$ and $\eta_4 \leq 2\eta_5 \leq 2\eta_6 \leq 2\eta_2$.
\end{theo}

Let us finish this section by showing that $\Gamma$ fulfils the all requirements for a genuine measure of noncompactness listed at the beginning. 

\begin{prop}\label{Konvex}
The function $\Gamma$ defined on bounded subsets of $X$ has the following properties:
\begin{enumerate}
\item $\Gamma(\overline{A})=\Gamma(A)$;
\item $\Gamma(A)=0$ if and only if $\overline{A}$ is SWC;
\item $\Gamma(A \cup B) = \max\{ \Gamma(A), \Gamma(B) \}$;
\item $\Gamma(\lambda A) = |\lambda| \Gamma(A)$ for all $\lambda \in {\Bbb R}$;
\item $\Gamma(A+B) \leq \Gamma(A) + \Gamma(B)$;
\item $\Gamma(\conv(A)) \leq 4 \, \Gamma(A)$.
\end{enumerate}
\end{prop}

\noindent
{\bf Proof.} (1) and (2) follow straightly from the definition of $\Gamma$. (3), (4) and (5) follow from set identities: 
$(A \cup B)^{\mathcal U} = A^{\mathcal U} \cup B^{\mathcal U}$, $(\lambda A)^{\mathcal U}= \lambda A^{\mathcal U}$ and 
$(A+B)^{\mathcal U} = A^{\mathcal U} + B^{\mathcal U}$. Statement (6), the most tricky, was proved in \cite[Theorem 6.7]{GR}.\findemo

\section{Quantifying uniform convexity for operators}

In this section we will discuss the application of the measure of weak noncompactness. For an operator 
$T: X \rightarrow Y$, we will write $\Gamma(T):=\Gamma(T(B_X))$. We have the following.

\begin{prop}\label{interpol}
Let $A \subset X$ be convex symmetric bounded set with $\Gamma(A) <\varepsilon$. Then there exists a Banach space $Z$ and an operator $T: Z \rightarrow X$ such that $\|T\|=1$, $A \subset T(B_Z)$ and $\Gamma(T) <\varepsilon$.
\end{prop}

\noindent
{\bf Proof.} Without loss of generality we may assume that $A$ is closed. Then, just take $Z = \spn (A)$, endow it with the norm given by the Minkowski functional of $A$ and take $T$ the identity operator.\findemo\\

If we consider the alternative measure of weak noncompactness $\Gamma_s$ introduced after Proposition \ref{superweaknoncom}, we have the following quantified version of the symmetry of the operator ideal ${\mathfrak W}^{super}$.

\begin{theo}\label{duality}
Let $T: X \rightarrow Y$ and operator. Then $\Gamma_s(T^*)=\Gamma_s(T)$.
\end{theo}

\noindent
{\bf Proof.} We will assume firstly that $ \Gamma_s(T)>0$.
Take $0<\varepsilon < \Gamma_s(T)$. Then, for every $N \in {\Bbb N}$ there are elements 
$(x_n)_{n=1}^N \subset B_X$ and $(x^*_n)_{n=1}^N \subset B_{X^*}$ such that 
$$\langle x^*_n , T(x_m) \rangle =0 \mbox{ for } m < n, $$
$$\langle x^*_n , T(x_m) \rangle \geq \varepsilon \mbox{ for } m \geq n .$$
But this is exactly the same that
$$\langle T^*(x^*_n) ,  x_m  \rangle =0 \mbox{ for } m < n, $$
$$\langle T^*(x^*_n) , x_m \rangle \geq \varepsilon \mbox{ for } m \geq n .$$
By reversing the order of $1,\dots, N$, we get $\Gamma_s (T^*) \geq \varepsilon$. That gives 
$\Gamma_s(T^*) \geq\Gamma_s(T)$. 
Suppose now that $ \Gamma_s(T^*)>0$ and take $0<\varepsilon < \Gamma_s(T^*)$.
Then, for every $N \in {\Bbb N}$ there are elements 
$(x^{**}_n)_{n=1}^N \subset B_{X^{**}}$ and $(x^*_n)_{n=1}^N \subset B_{X^*}$ such that 
$$\langle x^{**}_n , T^*(x^*_m) \rangle =0 \mbox{ for } m < n, $$
$$\langle x^{**}_n , T^*(x_m) \rangle \geq \varepsilon \mbox{ for } m \geq n .$$
Fix $\lambda>1$. Helly's theorem \cite[p. 159]{banach} allows us to find $(x_n)_{n=1}^N \subset \lambda B_X$ 
such that 
$$  \langle x^{**}_n , T^*(x^*_m) \rangle = \langle x_n , T^*(x^*_m)  \rangle $$
for every $1 \leq n,m \leq N$. That implies $\Gamma_s(T) \geq \lambda^{-1}\varepsilon$, after reversing the order of $1,\dots, N$. By the arbitrarily choice of constants, we get $\Gamma_s(T) \geq \Gamma_s(T^*)$.\\
So far we have proved that  $ \Gamma_s(T)>0$ if and only if  $ \Gamma_s(T^*)>0$ and, in such a case, 
$ \Gamma_s(T)= \Gamma_s(T^*)$. That also implies  $ \Gamma_s(T)=0$ if and only if  $ \Gamma_s(T^*)=0$, therefore the proof is complete.\findemo

\begin{coro}\label{duality2}
Let $T: X \rightarrow Y$ be an operator. Then $2^{-1}\Gamma(T) \leq \Gamma(T^*) \leq 2\Gamma(T)$.
\end{coro}

\begin{rema}
Using $\gamma_2$ as a measure of weak noncompactness for sets and operators, the quantified version of Gantmacher theorem \cite{AC} would become an equality.
\end{rema}

De Blasi's measure applied to operators does not satisfy a similar quantified Gantmacher result, as observed in \cite{AC} after an example from \cite{Asta}, neither does the measure on super weak noncompactness introduced by Tu \cite{KT2}, inspired by De Blasi's definition, as
$$ \sigma(T) = \inf\{ \varepsilon>0: \exists K \subset Y, K \mbox{ is SWC}, T(B_X) \subset K + \varepsilon B_Y \} $$
Indeed, Tu shows provides a sequence of operators $T_n$ such that  and $\sigma(T^*_n)=1$ for all $n \in {\Bbb N}$ and $\lim_n \sigma(T_n) =0$.\\

Now we will consider a notion of uniform convexity for operators.
In order to make notation shorter, for a convex function $f$ we will write
$$ \Delta_f(x,y) =  \frac{f(x)+f(y)}{2} - f \left( \frac{x+y}{2} \right) . $$
An operator $T: X \rightarrow Y$ is called uniformly convex if for every $\varepsilon>0$ there is a $\delta>0$ such that 
$\|T(x)-T(y)\| \leq \varepsilon$ whenever $x,y \in B_{X}$ are such that  $\Delta_{\| \cdot \|^2}(x,y) <\delta$. An operator $T: X \rightarrow Y$  is called uniformly convexifying if it becomes uniformly convex after a suitable renorming of $X$. It turns out that the class of uniformly convexifying operators agrees with ${\mathfrak W}^{super}$.\\

We will say that $T$ is $\varepsilon$-uniformly convex ($\varepsilon$-UC) if there is $\delta>0$ such that $\|T(x)-T(y)\| \leq \varepsilon$ whenever $x,y \in B_X$ are such that $\Delta_{\| \cdot \|^2} (x,y) < \delta$. The following result contains two alternative forms of the $\varepsilon$-UC property that we will need later.

\begin{lema}\label{epUC}
	For an operator  $T: X \rightarrow Y$ and $\varepsilon>0$, the following statements are equivalent:
	\begin{itemize}
		\item[(i)] $T$ is $\varepsilon$-UC;
		\item[(ii)] $\limsup_n \|T(x_n) - T(y_n)\| \leq \varepsilon$ whenever $x_n, y_n \in B_X$ are such $$\lim_n \Delta_{\| \cdot \|^2}(x_n,y_n) =0;$$ 
		\item[(iii)] there is $\delta>0$ such that $\|T(x)-T(y) \| \leq \varepsilon$ whenever $x,y \in X$ are such that $\|x\|=\|y\|=1$ and $\|x+y\|>2(1-\delta)$.
	\end{itemize}
\end{lema}

\noindent
{\bf Proof.} The proof is left to the reader.\findemo\\

For the construction of a quantified uniformly convex norm we will use this result.

\begin{theo}[\cite{GR}] \label{renorming}
	Let $(X,\| \cdot \|)$ be a Banach space,
	let $f: X \rightarrow [0,+\infty]$ be a proper convex function and let $C \subset \dom(f)$ be a bounded convex set.
	Assume $f$ is Lipschitz on $C$.
	Then given $\delta>0$ there exists an equivalent norm $|\!|\!|  \cdot |\!|\!| $ on $X$ and $\zeta>0$ such that 
	$\Delta_f(x,y) <\delta$ whenever $x,y \in C$ satisfy $\Delta_{ |\!|\!|  \cdot |\!|\!|^2 }(x,y) <\zeta$.
	Therefore, if $f$ 
	was moreover $\varepsilon$-uniformly convex for some $\varepsilon>0$ (with respect to a pseudo-metric) on $C$, then $|\!|\!|  \cdot |\!|\!|^2$ would be 
	$\varepsilon$-uniformly convex on $C$ (with respect to the same pseudo-metric).
\end{theo}

We are ready to prove the quantified Beauzamy's renorming result. 

\begin{theo}\label{eprenorm}
	Let $(X, \| \cdot \|)$ be a Banach space, and let $T:X \rightarrow Y$ be an operator such that $\Gamma(T)<\varepsilon$. Then there exists an equivalent norm $|\!|\!|  \cdot |\!|\!| $ on $X$ such that 
	$|\!|\!|  \cdot |\!|\!|  \leq \| \cdot \|$ and such that $T$ is $\varepsilon$-UC on $(X,|\!|\!|  \cdot |\!|\!| )$.\\ 
	Moreover, in case $X$ and $Y$ are dual Banach spaces and $T$ is an adjoint operator, then the norm 
	$|\!|\!|  \cdot |\!|\!| $ making  $T$ is $\varepsilon$-UC can be taken to be a dual one.
\end{theo}

\noindent
{\bf Proof.} Take $\Gamma(T) = \varepsilon' < \varepsilon$ 
and $1 < \lambda <\varepsilon/\varepsilon'$.
By Theorem  \ref{quantify-swc}, the set $B = \lambda \overline{T(B_X)}$ supports a convex bounded $\varepsilon$-uniformly convex function $f$ that we may assume it is also Lipschitz, see \cite[Proposition 5.4]{GR}. The function $f \circ T$ is $\varepsilon$-uniformly convex with respect to the pseudo-metric $d(x,y) =\| T(x) - T(y) \|$ on $\lambda B_X$. By Theorem \ref{renorming}, there is an equivalent norm  $\| \cdot \|_u$ on $X$ that is $\varepsilon$-uniformly convex with respect to $d$ on the set $\lambda B_X$. All the norms defined by the formula
$$  |\!|\!|  \cdot |\!|\!| ^2 = \lambda^{-2} \| \cdot \|^2 + \xi \| \cdot \|_u^2 $$
are $\varepsilon$-uniformly convex with respect to $d$ on the set $\lambda B_X$. By taking $\xi>0$ small enough we may assume that 
$$ \lambda^{-1} \| \cdot \| \leq |\!|\!|  \cdot |\!|\!|  \leq \| \cdot \|.$$
Since the unit ball of $ |\!|\!|  \cdot |\!|\!| $ contains  $\lambda B_X$, we get that 
$T$ becomes $\varepsilon$-UC  when $X$ is endowed with $ |\!|\!|  \cdot |\!|\!| $.\\
Assume now that $X$ and $Y$ are dual spaces and $T$ is an adjoint operator, and therefore it is weak$^*$ to weak$^*$ continuous.
By the first part, we may assume that $X$ is already endowed with a (non dual) norm such that $T$ is $\varepsilon$-UC. We claim that the norm $|\!|\!|.|\!|\!|$ on $X$ having $\overline{B_X}^{w^*}$ as the unit ball makes $T$ $\varepsilon$-UC too. By Lemma \ref{epUC} there is $\delta>0$ such that $x,y \in B_X$
and $\|x+y\| > 2(1-\delta)$ implies $\|T(x)-T(y)\| \leq \varepsilon$.
Therefore,  $\mbox{diam}(T(H \cap B_X)) \leq \varepsilon$ whenever $H$ is a halfspace such that
$H \cap (1-\delta)B_X=\emptyset$.
Take $x,y \in X$ with $|\!|\!|x|\!|\!|=|\!|\!|y|\!|\!|=1$ and $|\!|\!|x+y|\!|\!| > 2(1-\delta/2)$. Note that the condition implies that the segment $[x,y]$ does not meet $(1-\delta)\overline{B_X}^{w^*}$
Take $H$ a weak$^*$-open halfspace such that $[x,y] \cap (1-\delta)\overline{B_X}^{w^*}=\emptyset$. 
We have $\|x-y\| \leq  \mbox{diam}(H \cap \overline{B_X}^{w^*})$.
Now, by the weak$^*$ to weak$^*$-continuity of $T$ we have
\[ T(H \cap \overline{B_X}^{w^*}) \subset \overline{T(H \cap B_X)}^{w^*} . \]
As $\mbox{diam}(\overline{T(H \cap B_X)}^{w^*})=\mbox{diam}(T(H \cap B_X)) \leq \varepsilon$ by the weak$^*$ semicontinuity of the norm of $Y$ and the previous observation, we get that $\| T(x) - T(y)\| \leq \varepsilon$ as wished.\findemo

\section{Proof of the main result and consequences}

The norm of the Banach space $(X,\| \cdot \|)$ is said {\it uniformly Gateaux smooth} if 
for every $h \in X$
\[ \sup\{ \|x+th\|+\|x-th\|-2: x \in S_X\} =o(t) \mbox{~when~} t \rightarrow 0. \]
It is well known \cite[Theorem 6.7]{DGZ} that the norm on $X$ is uniformly Gâteaux smooth if and only if the dual norm on $X^*$ is weak$^*$ uniformly rotund (W$^*$UR), that is, weak$^*$-$\lim_n (x^*_n - y^*_n)=0$ whenever $x^*_n, y^*_n \in B_{X^*}$ are such that $\lim_n \Delta_{\| \cdot \|^2} (x^*_n, y^*_n) =0$.

\begin{lema}\label{lemarenorm}
	Let $A \subset X$ be a subset and let $\varepsilon>0$. Assume that $A=\bigcup_{k=1}^\infty A_k$ with $A_k$ bounded and $\Gamma(A_k) < \varepsilon$ for every $k \in {\Bbb N}$. Then, there exists an equivalent  norm $|\!|\!|  \cdot |\!|\!| $ on $X$ such that the dual norm on $X$ has the following property: whenever $(x^*_n), (y^*_n) \subset B_{X^*}$ are such that $\lim_n \Delta_{\| \cdot \|^2} (x^*_n, y^*_n) =0$, then
	$$ \limsup_n |x^*_n(x) - y^*_n(x)| \leq 8\varepsilon  $$
	for every $x \in A$.
\end{lema}

\noindent
{\bf Proof.} Let $B_k$ be the symmetric convex hull of $A_k$. By Proposition \ref{Konvex}, we have $\Gamma(B_k) <4\varepsilon$. Let $T_k: Z_k \rightarrow X$ the operator given by Proposition \ref{interpol} such that $\Gamma(T_k) <4\varepsilon$ and $A_k \subset B_k \subset T_k(B_{Z_k})$. Now, by Corollary \ref{duality2} $\Gamma(T_k^*) < 8\varepsilon$, and, by Theorem \ref{eprenorm}, $T_k^*$ became $8\varepsilon$-UC with an equivalent dual norm $\| \cdot  \|_k \leq \| \cdot \|$. Consider the equivalent dual norm on $X^*$ defined by the formula
$$ |\!|\!|  \cdot |\!|\!| ^2 = \sum_{k=1}^\infty 2^{-k} \| \cdot \|_k^2 .$$
Suppose given $(x^*_n), (y^*_n) \subset B_{X^*}$ with $\lim_n \Delta_{\| \cdot \|^2} (x^*_n, y^*_n) =0$.
Then, for every $k \in {\Bbb N}$, we have  $\lim_n \Delta_{\| \cdot \|^2_k} (x^*_n, y^*_n) =0$ and therefore $\limsup_n \| T_k^*(x^*_n) - T_k^*(y^*_n) \| \leq 8\varepsilon$ on $Z^*_k$. In particular, for every $z \in Z_k$, we get 
$$  \limsup_n | \langle T_k(z), x^*_n \rangle - \langle T_k(z), y^*_n \rangle | = 
\limsup_n | \langle z, T_k^*(x^*_n) \rangle - \langle z, T_k^*(y^*_n) \rangle | \leq 8\varepsilon. $$
Having in mind that $A_k \subset T(B_{Z_k})$, we obtain 
$  \limsup_n | x^*_n(x)  -  y^*_n(x) | \leq 8\varepsilon $
for every $x \in A_k$. Since this is true for every $k \in {\Bbb N}$, the lemma is proved.
\findemo\\

Instead of proving Theorem \ref{main}, we will prove of the following, being both equivalent thanks to \cite{BRW}.

\begin{theo}\label{main2}
	Let $X$ be a Banach space. The following statements are equivalent:
	\begin{itemize}
		\item[(i)] $X$ is a subspace of a Hilbert generated space;
		\item[(ii)] For every $\varepsilon>0$ there are sets $(B_n^\varepsilon)$ such that $B_X=\bigcup_{n=1}^\infty B_n^\varepsilon$ and $\Gamma(B_n^\varepsilon) < \varepsilon$;
		\item[(iii)] There exists a linearly dense set $A \subset X$ such that for every $\varepsilon>0$ it can be decomposed as $A=\bigcup_{n=1}^\infty A_n^\varepsilon$ where each $A_n^\varepsilon$ is bounded and $\Gamma(A_n^\varepsilon) < \varepsilon$;
		\item [(iv)] $X$ admits an equivalent uniformly Gâteaux norm.
	\end{itemize}
\end{theo}

\noindent
{\bf Proof.} (i)$\Rightarrow$(ii) It is enough to prove statement (ii) for a Hilbert generated space since that property is clearly inherited by subspaces. Let $H$ be a Hilbert space and $T: H \rightarrow X$ an operator with dense range. For every $0 < \varepsilon' < \varepsilon$ we have 
$$ B_X \subset \bigcup_{n=1}^\infty (nT(B_H) + \varepsilon' B_X ). $$
We have $\Gamma(nT(B_H) + \varepsilon' B_X ) \leq \varepsilon'$ and we can take $B_n^\varepsilon =B_X \cap
(nT(B_H) + \varepsilon' B_X )$.\\
(ii)$\Rightarrow$(iii) It is obvious.\\
(iii)$\Rightarrow$(iv) By Lemma \ref{lemarenorm}, for every $k \in {\Bbb N}$ there exists an equivalent dual norm $\| \cdot \|_k$ on $X^*$ such that: whenever $(x^*_n), (y^*_n) \subset B_{X^*}$ are such that $\lim_n \Delta_{\| \cdot \|_k^2} (x^*_n, y^*_n) =0$, then
$$ \limsup_n |x^*_n(x) - y^*_n(x)| \leq 1/k $$
for every $x \in A$. The dual norm defined by
$$ |\!|\!|  \cdot |\!|\!| ^2 = \sum_{k=1}^\infty 2^{-k}\| \cdot \|_k^2 $$
satisfies then $ \limsup_n |x^*_n(x) - y^*_n(x)| =0 $ whenever $x \in \spn(A)$ and  $(x^*_n), (y^*_n) \subset B_{X^*}$ are such that $\lim_n \Delta_{ |\!|\!|  \cdot |\!|\!| ^2} (x^*_n, y^*_n) =0$. As the sequences  $(x^*_n), (y^*_n)$ are bounded and $\spn(A)$ is dense, we have  $ \limsup_n |x^*_n(x) - y^*_n(x)| =0 $ for every $x \in X$. Therefore, the norm $ |\!|\!|  \cdot |\!|\!| $ is W$^*$UR and its predual norm on $X$ is uniformly Gâteaux.\\
(iv)$\Leftrightarrow$(i) It was proved in \cite{FGZ} (see also \cite[Theorem 6.30]{BIO}).\findemo\\

The result of Fabian, Godefroy and Zizler \cite{FGZ}  (see also \cite[Theorem 6.30]{BIO}) gives actually more information: the linearly dense set can be decomposed, for every $\varepsilon >0$ in countably many pieces which are uniformly weakly null up to $\varepsilon$ in the sense of Proposition \ref{Banach-Saks}. That cannot be done on every set generally, however it applies to Markushevich bases as a consequence of the following dual interpretation of \cite{farma}.

\begin{prop}\label{AFnull}
Let $X$ be a subspace of a Hilbert generated Banach space and let $A \subset X$ be a bounded set such that $0$ is its only cluster point and $A \cup \{0\}$ is weakly compact. Then, for every $\varepsilon>0$ there is a decomposition  $A=\bigcup_{n=1}^\infty A_n^\varepsilon$ such that for every $n \in {\Bbb N}$ and for every $x^* \in B_{X^*}$ then
$$ | \{ x \in A_n^\varepsilon: |x^*(x)| > \varepsilon \} | \leq n .$$
\end{prop}

\noindent
{\bf Proof.} Without loss of generality we may assume $A \subset C(K)$ where $K$ is uniform Eberlein.
Indeed, take $K=(B_{X^*},w^*)$ that is uniform Eberlein after a result from \cite{BRW}. Consider the embedding of $K$ into $\ell_\infty(A)$ given by $K \ni t \rightarrow (f(t))_{f \in A}$ and note that, actually, it take values into $c_0(A)$. By \cite{farma} (see also \cite[Theorem 6.33]{BIO}), there is a decomposition of the index set  $A=\bigcup_{n=1}^\infty A_n^\varepsilon$ such that for every 
$t \in K$
$$ | \{ f \in A_n^\varepsilon: |f(t)| > \varepsilon \} | \leq n .$$
Since $K$ is a norming set on $C(K)$, we get the conclusion for every norm one functional, see Remark \ref{norming}.\findemo\\

It is interesting to investigate the case where statement (iii) Theorem \ref{main2} happens without countable decomposition of the linearly dense set, that is, when that set is relatively SWC. Let us recall that a Banach space $X$ that contains a linearly dense SWC set is called {\it super weakly compactly generated} (super WCG). That condition is equivalent to the existence of a super weakly compact operator into $X$ with dense range. See \cite{raja2} for the renorming properties of super WCG Banach spaces.

\begin{theo}
	Let $X$ be a Banach space. The following are equivalent:
	\begin{itemize}
		\item[(i)] $X$ is super WCG;
		\item[(ii)] $X$ has a  Markushevich basis $\{x_i,x^*_i\}_{i \in I}$ such that $\{x_i:i \in I\} \cup \{0\}$ is SWC;
		\item[(iii)] There exists an one-to-one bounded linear operator $T: X^* \rightarrow c_0(I)$, for some set $I$, which is weak$^*$ to pointwise continuous and SWC.
	\end{itemize}	
\end{theo}

\noindent
{\bf Proof.} (i)$\Rightarrow$(ii) Without loss of generality we may assume that $X$ is generated by a balanced convex SWC set $K$. The proof of existence of Markushevich basis on WCG spaces allows the choice  $\{x_i:i\in I\} \subset K$, see \cite[Theorem 13.16]{banach}. Clearly, the only cluster point of  $\{x_i:i \in I\}$ is $0$, and thus  $\{x_i:i \in I\} \cup \{0\}$ is SWC.\\
(ii)$\Rightarrow$(iii) We may assume that  $\{x^*_i:i \in I\}$ is uniformly bounded. Define $T(x^*)=(x^*(x_i))_{i \in I}$ which, initially, takes values into $\ell_\infty(I)$. It can be proved that $T(X^*) \subset c_0(I)$, see \cite[Theorem 12.20]{banach} for the details. In order to see that $T$ is SWC, we will see that $T^*$ is SWC. Indeed, $T^*$ takes the basis $(e_i)_{i\in I}$ of $\ell_1(I)$ to the set $\{x_i:i\in I\}$. Since  $B_{\ell_1(I)}$ is the closed convex hull of $(e_i)_{i\in I}$, we deduce that 
$T(B_{\ell_1(I)})$ is contained in the balanced convex hull of  $\{x_i:i\in I\}$, and therefore it is relatively SWC.\\ 
(iii)$\Rightarrow$(i) Consider the adjoint operator $T^{*}: \ell_1(I) \rightarrow X^{**}$, which is SWC, and note that every element from $(e_i)_{i\in I}$, the basis of $\ell_1(I)$, goes through $T^{**}$ to a weak$^*$ continuous element of $X^{**}$. Therefore $T^{*}(\{e_i:i\in I\}) \subset X$, and thus $T^{*}(\ell_1(I)) \subset X$. Now, as $T$ is one-to-one, $T^*$ has a dense range and therefore $X$ is super WCG.\findemo

\begin{rema}
Recall that being $X$ super WCG is equivalent to the existence of an equivalent strongly uniformly Gâteaux norm on $X$ \cite[Theorem 1.6]{raja2}.
\end{rema}

Now we will give an application to {\it Jordan algebras}. We refer the reader to \cite{HKPP} for the necessary definitions.
In \cite{HKPP} the authors have proved that the measures of weak noncompactness $\gamma$ and $\omega$ (De Blasi's measure) agree on a {\it JBW$^*$-triple predual}. The next result shows that we can add $\Gamma$ to them.

\begin{prop}
Let $X$ be JBW$^*$-triple predual. Then $\omega$, $\gamma$ and $\Gamma$ agree on $X$.
\end{prop}

\noindent
{\bf Proof.} Let $A \subset X$ be bounded and take $\varepsilon>\gamma(A)$. Since $\gamma=\omega$ by \cite{HKPP}, there is $K \subset X$ weakly compact such that $A \subset K + \varepsilon B_X$. By \cite[Theorem 6.3]{LR}, $K$ is SWC. Therefore, the inclusion $A \subset K + \varepsilon B_X$ implies $\Gamma(A) \leq \varepsilon$. We deduce $\Gamma(A) \leq \gamma(A)$. Since the other inequality always holds $\Gamma(A) = \gamma(A)=\omega(A)$.\findemo\\

This result implies for a  JBW$^*$-triple predual that the notions of WCG and super WCG are equivalent. Moreover, in \cite[Theorem 9.3]{HKPP} the authors provide characterizations for JBW$^*$-triple predual to be WCG or strongly WCG. It turns out that in such cases the spaces become super WCG or strongly super WCG (S$^2$WCG), respectively, which implies nice geometrical properties under renorming, see Theorem 1.6 and Theorem 1.9 in \cite{raja2}.

\section{More on uniformly weakly null sets}

A set that contains a sequence equivalent to the basis of $\ell_1$ cannot be uniformly weakly null. Whether a Schauder basis is a uniformly weakly null set or not will be characterized among the symmetric basis. Recall that an unconditional Schauder basis is said {\it symmetric} if it is uniformly equivalent to all its permutations.
The following is a result  due to Troyanski \cite{tro1} reformulated in our terms.

\begin{theo}
	Let $X$ be a Banach space with a symmetric basis $(e_i)_{i\in I}$. Then the following statements are equivalent:
	\begin{itemize}
		\item[(i)]  $\{e_i : i \in I\}$ is uniformly weakly null;
		\item[(ii)] $\{e_i : i \in I\} \cup \{0\}$ is SWC;
		\item[(iii)] $0$ is a weak cluster point of $\{e_i : i \in I\}$;
		\item[(iv)]$X$ is not isomorphic to $\ell_1(I)$.
	\end{itemize}
In case $I$ is not countable (equivalently, $X$ is not separable), these conditions characterize the existence of an equivalent uniformly Gâteaux norm on $X$.
\end{theo}

\noindent
{\bf Proof.} Note that (i)$\Leftrightarrow$(iv) and characterization of uniform Gâteaux renorming for Banach spaces with symmetric bases is the original result of Troyanski \cite{tro1}, see also \cite[Lemma 7.52]{BIO} and  \cite[Theorem 7.54]{BIO}.
Clearly (i)$\Rightarrow$(ii) and (i)$\Rightarrow$(iii). On the other hand, assume (iii) and let $c_s \geq 1$ the symmetric unconditionality constant of the basis. For every $\varepsilon>0$ there are indices $(i_k)_{k=1}^n \subset I$ and positive numbers $\lambda_k$, $1\leq k \leq n$, with $\sum_{k=1}^n \lambda_k=1$ such that
$$ \| \sum_{k=1}^n \lambda_k e_{i_k} \| \leq \varepsilon .$$
Consider a cyclic permutation $\lambda_k^1=\lambda_{k+1}$ if $k<n$ and $\lambda_n^1=\lambda_1$. Then we have
$$ \| \sum_{k=1}^n \lambda_k^1 e_{i_k} \| \leq c_s \varepsilon .$$
If we take the other $n-2$ cyclic permutations obtained by iterating the first one, the sum gives 
$$  \| \sum_{k=1}^n  e_{i_k} \| \leq  \| \sum_{k=1}^n \lambda_k e_{i_k} \| +
 \| \sum_{k=1}^n \lambda_k^1 e_{i_k} \| + \dots +  \| \sum_{k=1}^n \lambda_k^{(n-1)} e_{i_k} \|
\leq  n \, c_s \, \varepsilon . $$
Again, by the symmetry of the basis, for any $J \subset I$ with $n$ elements we have
$$  n^{-1} \| \sum_{i \in J}  e_{i} \| \leq c_s^2 \, \varepsilon $$
that implies $\{e_i : i \in I\}$ is uniformly weakly null. Finally, assume (ii). By Corollary \ref{Mercou2}, there exists an infinite sequence in the set  $\{e_i : i \in I\}$ which is uniformly weakly convergent. Since the unique allowed cluster point is $0$, the sequence is uniformly weakly null. That behavior can easily be extended to all the basis  $(e_i)_{i\in I}$ by the symmetry.\findemo\\

Without the hypothesis of symmetry for the basis, we have the following result.

\begin{prop}
Let $X$ be a Banach space with nontrivial type. Then every unconditional seminormalized basic sequence (or set) is uniformly weakly null.
\end{prop}

\noindent
{\bf Proof.}  Let $(e_n)$ be an unconditional basic sequence with unconditionality constant $c_u \geq 1$, let $p \in (1,2]$ be the type of $X$ and $c_\tau$ the type constant. Without loss of generality we may assume $(e_n)$ is normalized. We have
$$ \| \sum_{n \in F} e_n \|  \leq c_u \, \| \sum_{n \in F} \epsilon_n e_n \| $$
whenever $\epsilon_n \in \{-1,1\}$ and $F \subset {\Bbb N}$ finite. Let $(r_n(t))$ denote the sequence of Rademacher functions. Applying the definition of type we get
$$ \| \sum_{n \in F} e_n \|  \leq c_u \, \int_0^1 \| \sum_{n \in F} r_n(t) \, e_n \| \, dt 
\leq c_u c_\tau  ( \sum_{n \in F} \| e_n \|^p )^{1/p} = c_u c_\tau \, n^{1/p} $$
that implies $(e_n)$ is a uniformly weakly null set.\findemo\\

In \cite{LR} it is proved a result about the coordinate combinatoric behavior of the SWC compact subsets of $c_0({\Bbb N})$ that are made up of characteristic functions. The following result shows that uniformly weakly null subsets in $c_0(I)$ made up of characteristic functions are more boring.

\begin{prop}
Let $\mathcal{F}$ be a family of finite subsets of a set $I$. 
Then $A=\{ \chi_F: F \in \mathcal{F}\}$ is uniformly 
weakly null as a subset of $c_0(I)$ if and only if there is a Hilbert space $H$ and an operator $T: H \rightarrow c_0(I)$ such that $A$ is covered by the image of an orthonormal basis of $H$. Moreover, an analogous result fails if $c_0(I)$ is replaced by another space with a long unconditional basis.
\end{prop}

\noindent
{\bf Proof.}
Consider the Hilbert space $H=\ell_2(A)$ with the basis $\{e_x: x \in A\}$.
As $A$ is uniformly weakly null, there exists $N$ such that 
$$ |\{i \in I: |x_i| >0 \}| \leq N $$
for every $x=(x_i)_{i \in I} \in A$. That implies that the assignation $e_x \rightarrow x$ can be extended to a linear operator. Indeed, for $(a_x)_{x \in A} \subset {\Bbb R}$ finitely supported, the sum $\sum_{x \in A} a_x x$ takes values in $c_0(I)$ and the bound  
$$ | (\sum_{x \in A} a_x x)_{i} | \leq \sum_{x_i \not =0} |a_x|  \leq N \sup_{x \in A}\{|a_x|\} 
\leq N \| \sum_{x \in A} a_x e_x \|$$
implies that the operator 
$$ T ( (a_x)_{x \in A}) = \sum_{x \in A} a_x x$$ 
can be extended to all $(a_x)_{x \in A} \in H$
with $\|T\| \leq N$. On the other hand, if the set $A$ is covered by the image of an orthonormal basis of a Hilbert space $H$ through an operator $T$, fix for every $x \in A$ an element $e_x \in H$ such that $(e_x)_{x \in A}$ is orthonormal. The fact that $A$ is uniformly weakly null follows easily from the fact that 
$$ \|x_1+ \dots + x_n \| \leq \|T\| \, \| e_{x_1} + \dots + e_{x_n} \| = \|T\| \, n^{1/2} $$
for different points $x_1,\dots,x_n \in A$.
For the last statement, we claim that $c_0(I)$ cannot be replaced by $\ell_{3/2}(I)$. Indeed, the canonical basis of $\ell_{3/2}(I)$ is a uniformly weakly null set, however it cannot be covered by the image of an operator from a Hilbert space (for $I$ uncountable). Otherwise $\ell_{3/2}(I)$ would be Hilbert generated, which is not the case by \cite[p. 316]{FGHZ}.\findemo

\begin{rema}
According to a classic result of Davis, Johnson, Lindenstrauss and  Pełczyński \cite{DJLP}, every relatively weakly compact set whose unique accumulation point is $0$ (like as in the hypothesis of Proposition \ref{AFnull}) is the image  through an operator of an unconditional basis in a reflexive space.
\end{rema}

The second named author proved that convex SWC sets considered with the weak topology are uniformly Eberlein \cite{raja}. These last results will deal with a more restrictive property. Following \cite{FGZ},
we say that a compact subset $K \subset E$ in a locally convex space is {\it linearly uniformly Eberlein} if there exists a linear injection $T:E \rightarrow c_0(I)$ which is continuous to the pointwise topology of $c_0(I)$ and for every $\varepsilon>0$ there is $n(\varepsilon)$  such that 
$$ |\{ i \in I: |T(x)_i| > \varepsilon \} | \leq n(\varepsilon) $$
for every $x \in K$. In case $K$ is moreover convex, we say that $K$ is affinely uniformly Eberlein if an affine map can be defined on $K$ with values on $c_0(I)$ with similar properties.

\begin{prop}
Let $X$ be a Banach space. Then $X$ contains a linearly dense uniformly weakly null set if and only if $(B_{X^*},w^*)$ is linearly uniformly Eberlein.
\end{prop}

\noindent
{\bf Proof.}
Let $A \subset X$ be uniformly weakly null. Observe that $T:X^* \rightarrow c_0(A)$ given by $T(x^*)=(x^*(x))_{x\in A}$ is well defined, one-to-one, and linearly represents $B_{X^*}$ as uniform Eberlein.
On the other hand, if $T: X^* \rightarrow c_0(I)$ witnesses that $(B_{X^*},w^*)$ is linearly uniformly Eberlein, then the coordinate maps define elements $\{x_i : i \in I \} \subset X$, as they are weak$^*$ continuous. It is not difficult to check that  $\{x_i : i \in I \}$ is uniformly weakly null and linearly dense.\findemo

\begin{rema}
Note that $(B_{X^*},w^*)$ can be uniformly Eberlein but not  linearly uniformly Eberlein. For that, just take a uniformly Gâteaux Banach space which is not WCG, for instance, Rosenthal's non WCG subspace of some  $L_1(\mu)$ space, see \cite{FGHZ} for more details. 
\end{rema}

This is the main question we cannot answer with the techniques of this paper.

\begin{open}
Is every super WCG Banach space generated by a uniformly weakly null set?
\end{open} 

Next result is motivated by  \cite[Theorem 4]{FGHZ}.

\begin{theo}
	Let $K \subset X$ be a super weakly compact convex subset of density (equivalently, weight) $\omega_1$. Then $K$ is affinely  uniformly Eberlein. 
\end{theo}

\noindent
{\bf Proof.} Let $Z$ be a reflexive Banach space and $T:Z \rightarrow X$ a one-to-one super weakly compact operator such that $K \subset T(B_Z)$, \cite[Theorem 1.3]{raja}. Then $K$ is linearly homeomorphic to a weakly compact subset of $Z$ of density $\omega_1$. Without loss of generality we may assume that $Z$ has density $\omega_1$ too.
Now, $T^*: X^* \rightarrow Z^*$ is a super weakly compact  operator with dense range. We deduce that $Z^*$ is super WCG and thus it is a uniformly Gâteaux renormable Banach space of density $\omega_1$. By \cite[Theorem 4]{FGHZ}, $B_Z$ is linearly uniformly Eberlein, which implies that $K$ is affinely uniformly Eberlein.\findemo\\

There are stronger results for weakly compact convex sets of weight strictly less than $\omega_1$, that is, the compact is metrizable. For instance, Keller's theorem, see \cite{BP}, showing an affine homeomorphism to the Hilbert cube.


\begin{thebibliography}{99}

\bibitem{AC}{\sc C. Angosto, B. Cascales},
{Measures of weak noncompactness in Banach spaces},
{\it Topology Appl.} 156 (2009), no. 7, 1412--1421. 

\bibitem{ADL}{\sc J.M. Ayerbe Toledano, T. Domínguez Benavides, G. López Acedo},
{\it Measures of noncompactness in metric fixed point theory}, 
Operator Theory: Advances and Applications, 99. Birkhuser Verlag, Basel, 1997.

\bibitem{farma}{\sc S. Argyros, V. Farmaki},
	{On the structure of weakly compact subsets of Hilbert spaces and applications to the geometry of Banach spaces}, {\it Trans. Amer. Math. Soc.} 289 (1985), 409--427.


\bibitem{Asta}{\sc K. Astala}, 
{On Measures of Noncompactness and Ideal Variations in Banach Spaces}, 
{\it Ann. Acad. Sci. Fenn. Math. Diss.} 29 (1980), 42 pages.

\bibitem{beau1} {\sc B. Beauzamy}, 
{Op\'erateurs uniform\'ement convexifiants},
{\it Studia Math.} {57} (1976), no. 2, 103--139.

\bibitem{beau2} {\sc B. Beauzamy},
{Quelques propri\'et\'es des op\'erateurs uniform\'ement convexifiants},
{\it Studia Math.} 60 (1977), no. 2, 211--222. 


\bibitem{beau}{\sc B. Beauzamy}
{\it Introduction to Banach spaces and their geometry},
North-Holland Mathematics Studies, 68.
Notas de Matem\'atica [Mathematical Notes], 86. North-Holland Publishing Co.,
Amsterdam-New York, 1982.

\bibitem{BRW}{\sc Y. Benyamini, M.E. Rudin, M. Wage},
{Continuous images of weakly compact subsets of Banach spaces},
{\it Pacific J. Math.} 70 (1977), no. 2, 309--324. 

\bibitem{BS}{\sc Y. Benyamini, T. Starbird},
{Embedding weakly compact sets into Hilbert space},
{\it Israel J. Math.} 23 (1976), no. 2, 137--141. 

\bibitem{BP}{\sc C. Bessaga, A. Pelczynski,}{\it ~Selected topics in
infinite-dimensional topology}, Mono. Mat. Tom {58}, PWN-Polish
Scientific Publishers, Warszawa, 1975.


\bibitem{CMR}{\sc B. Cascales, W. Marciszewski, M. Raja},
{Distance to spaces of continuous functions}. 
{\it Topology Appl.} 153 (2006), no. 13, 2303--2319. 


\bibitem{CPR}{\sc
B. Cascales, A. Pérez, M. Raja},
Radon--Nikodým indexes and measures of weak noncompactness. 
{\it J. Funct. Anal.} 267 (2014), no. 10, 3830--3858. 

\bibitem{CauseyDilworth} {\sc R.M. Causey, S. Dilworth},
{Metric characterizations of super weakly compact operators},
 Studia Math. 239 (2017), no. 2, 175--188.

\bibitem{Cheng}{\sc L. Cheng, Q. Cheng, B. Wang, W. Zhang},
{On super-weakly compact sets and uniformly convexifiable sets},  
{\it Studia Math.} 199 (2010), no. 2, 145--169. 


\bibitem{Cheng2}{\sc L. Cheng, Q. Cheng, J. Zhang},
{On super fixed point and super weak compactness of convex subsets in Banach spaces},
{\it J. Math. Anal. Appl.} 428 (2015), 1209--1224.

\bibitem{Cheng3}{\sc L. Cheng, Q. Cheng, S. Luo, K. Tu, J. Zhang},
{On super weak compactness of subsets and its equivalences in Banach spaces}
{\it J. Convex Anal.} 25 (2018), no. 3, 899--926.

\bibitem{DJLP}{\sc W.J. Davis, T. Figiel, W.B. Johnson, A. Pełczyński},
{Factoring weakly compact operators},
{\it J. Functional Analysis} 17 (1974), 311--327. 


\bibitem{DGZ}{\sc R. Deville, G. Godefroy, V. Zizler}
{\it Smoothness and renormings in Banach spaces},
Pitman Monographs and Surveys in Pure and Applied Mathematics, 64. Longman Scientific \& Technical,
Harlow, 1993.


\bibitem{DJP}{\sc J. Diestel, H. Jarchow, A. Pietsch},
{Operator Ideals},
{\it Handbook of the Geometry of Banach spaces}
Vol. 1, W.B. Johnson and J. Lindenstrauss editors,
Elsevier, Amsterdam (2001), 437--496.


\bibitem{Enflo}{\sc P. Enflo},
{Banach spaces which can be given an equivalent uniformly convex norm},
{\it Israel J. Math.} 13 (1972), 281--288.

\bibitem{ErMa}{\sc P. Erdős, M. Magidor},
{A note on regular methods of summability and the Banach-Saks property},
{\it Proc. Amer. Math. Soc.} 59 (1976), no. 2, 232--234. 

	
\bibitem{FGHZ}{\sc M. Fabian, G. Godefroy, P. H\'ajek, V. Zizler},
{Hilbert-generated spaces}, 
{\it J. Funct. Anal.} 200 (2003), no. 2, 301--323. 

\bibitem{FGZ}{\sc M. Fabian, G. Godefroy, V. Zizler},
{The structure of uniformly G\^ateaux smooth Banach spaces},
{\it Israel J. Math.} 124 (2001), 243--252. 


\bibitem{banach} {\sc M. Fabian, P. Habala, P. H\'ajek, V. Montesinos and V. Zizler},
{\it Banach Space Theory. The Basis for Linear and Nonlinear Analysis},
CMS Books in Mathematics,
Springer, New York, 2011.

\bibitem{FHMZ}{\sc M. Fabian, P. Hájek, V. Montesinos, V. Zizler}, 
A quantitative version of Krein's theorem, 
{\it Rev. Mat. Iberoamericana 21} (2005) 237--248.

\bibitem{FHMZ2}{\sc M. Fabian, P. Hájek, V. Montesinos, V. Zizler}, 
{Weakly compact generating and shrinking Markuševič bases},
{\it Serdica Math.} J. 32 (2006), no. 4, 277--288. 

\bibitem{FHZ}{\sc M. Fabian, P. H\'ajek, V. Zizler},
{Uniform Eberlein compacta and uniformly Gâteaux smooth norms},
{\it Serdica Math.} J. 23 (1997), no. 3-4, 351--362. 

\bibitem{FMZ}{\sc M. Fabian, P. Hájek, V. Montesinos, V. Zizler}, 
{A characterization of subspaces of weakly compactly generated Banach spaces},
{\it J. London Math. Soc.} (2) 69 (2004), no. 2, 457--464. 


\bibitem{granero}{\sc A. S. Granero}, 
An extension of the Krein-Šmulian theorem, 
{\it Rev. Mat. Iberoam.} 22(1) (2006), 93--110.


\bibitem{GR}{\sc G. Grelier, M. Raja},
{Uniformly convex functions},
{\it J. Math. Anal. Appl.} 505 (2022), Issue 1, 125442, 25 pp. 

\bibitem{HKPP}{\sc J. Hamhalter, O.F.K. Kalenda, A.M. Peralta, H. Pfitzner},
{Measures of weak non-compactness in preduals of von Neumann algebras and JBW$^*$-triples}, 
{\it J. Funct. Anal.} 278 (2020), no. 1, 108300, 69 pp. 

\bibitem{heinrich}{\sc S. Heinrich},
{Finite representability and super-ideals of operators},
{\it Dissertationes Math. (Rozprawy Mat.)} 172 (1980).

\bibitem{heinrich2}{\sc S. Heinrich},
{Ultraproducts in Banach space theory},
{\it J. Reine Angew. Math.} 313 (1980), 72--104. 

\bibitem{BIO}{\sc P. H\'ajek, V. Montesinos, J. Vanderwerff, V. Zizler},
{\it Biorthogonal systems in Banach spaces},
CMS Books in Mathematics/Ouvrages de Mathématiques de la SMC, 26. Springer, New York, 2008.

\bibitem{James1}{\sc R. C. James},
{Super-reflexive Banach spaces},
{\it Canad. J. Math.} {24} (1972), 896 -- 904.

\bibitem{JL}{\sc W. B. Johnson, J. Lindenstrauss},
{Basic concepts in the geometry of Banach spaces},
{\it Handbook of the Geometry of Banach spaces}
Vol. 1, W.B. Johnson and J. Lindenstrauss editors,
Elsevier, Amsterdam (2001), 1-- 84.

\bibitem{LPR}{\sc G. Lancien, A. Prochazka, M. Raja},
{Szlenk index of convex hulls}
{\it J. Funct. Anal.} 272 (2017), 498--521.

\bibitem{LR}{\sc G. Lancien, M. Raja},
{Non linear properties of super weakly compact sets},
{\it to appear Annal. Int. Fourier}

\bibitem{saks}{\sc J. Lopez-Abad, C. Ruiz, P. Tradacete},
{The convex hull of a Banach-Saks set},
{\it J. Funct. Anal.} 266 (2014), no. 4, 2251--2280. 

\bibitem{merco}{\sc S. Mercourakis},
{On Cesàro summable sequences of continuous functions}, 
{\it Mathematika} 42 (1995), no. 1, 87--104. 


\bibitem{Pfitzner} {\sc H. Pfitzner},
A conjecture of Godefroy concerning James' theorem. (English summary) 
{\it Q. J. Math.} 64 (2013), no. 2, 547--553. 

\bibitem{Pisier}{\sc G. Pisier},
{Martingales with values in uniformly convex spaces},
{\it Israel J. Math.} 20 (1975), no. 3-4, 326--350.

\bibitem{raja}{\sc M. Raja},
{Finitely dentable functions, operators and sets},
{\it J. Convex Anal.} {15} (2008), 219--233.


\bibitem{raja2}{\sc M. Raja},
{Super WCG Banach spaces},
{\it J. Math. Anal. Appl.} 439 (2016), no. 1, 183--196.

\bibitem{tro1}{\sc S. L. Troyanski},
{On nonseparable Banach spaces with a symmetric basis},
{\it Studia Math.} 53 (1975), 253--263. 

\bibitem{tro2}{\sc S. L. Troyanski},
{Uniform convexity and smoothness in every direction in nonseparable Banach spaces with unconditional bases},
{\it C. R. Acad. Bulgare Sci.} 30 (1977), no. 9, 1243--1246. 

\bibitem{KT}{\sc K. Tu}, 
{Convexification of super weakly compact sets and measure of super weak noncompactness},
{\it Proc. Amer. Math. Soc.} 149 (2021), no. 6, 2531--2538. 

\bibitem{KT2}{\sc K. Tu},
{Equivalence of semi-norms related to super weakly compact operators},
{\it Bull. Aust. Math. Soc.} 104 (2021), no. 3, 506--518. 

\bibitem{wenzel}{\sc J. Wenzel},
{Uniformly convex operators and martingale type},
{\it Rev. Mat. Iberoamericana} 18 (2002), no. 1, 211--230. 

\bibitem{YLC}{\sc Z.T. Yang, Y.F. Lu, Q.J. Cheng},
{Super weak compactness and uniform Eberlein compacta},
{\it Acta Math. Sin.} (Engl. Ser.) 33 (2017), no. 4, 545--553. 



\end{thebibliography}
\end{document}